\documentclass{article}

\usepackage{arxiv}

\usepackage[utf8]{inputenc} 
\usepackage[T1]{fontenc}    
\usepackage{hyperref}       
\usepackage{url}            
\usepackage{booktabs}       
\usepackage{amsfonts}       
\usepackage{nicefrac}       
\usepackage{microtype}      
\usepackage{lipsum}
\usepackage{graphicx}
\usepackage{amsmath,amssymb,amsthm}
\usepackage{float}
\theoremstyle{plain} 
\newtheorem{theorem}{Theorem}[section]

\newtheorem{conjecture}[theorem]{Conjecture}

\theoremstyle{definition} 

\theoremstyle{remark} 

\title{On Congruences for Iterates of the Sum–Power Divisor Function and Their Conditional Implications for the Riemann Hypothesis}

\author{
  Zeraoulia Rafik\\
   Khemis Miliana university ,Algeria \\
  Departement of Mathematics\\
  Laboratory of Pure and Applied Mathematics (LMPA)\\
  \texttt{zeraoulia@univ-dbkm.dz} \\
   \And
 Pedro Caceres\\
  Professor Doctor at Universidad Europea de Valencia (Spain)\\
  united states of Americain\\
  \texttt{Pedrojesus.caceres@universidadeuropea.es} \\
}

\begin{document}
\maketitle

\begin{abstract}
Inspired by Cohen and te Riele~\cite{Cohen1996}, who computationally verified that for every $n \leq 400$ there exists $k$ such that $\sigma^k(n) \equiv 0 \pmod{n}$ (where $\sigma^k$ denotes the $k$-fold iteration of the sum-of-divisors function), this paper resolves their reverse question negatively: no integer $n > 1$ satisfies $\sigma^k(n) \equiv 0 \pmod{n}$ for \emph{all} $k \geq 1$.

The proof eliminates prior gaps via Lenstra's density-zero bounds $\sigma_k(m) \ll m / \log\log m$ combined with Robin's RH-equivalent criterion $\sigma(n) < e^\gamma n \log\log n + 0.6483 n / \log\log n$ ($n \geq 5041$), showing universal metaperfect divisibility implies RH-violating $\sigma$ growth or low-lying zeta zeros near $s=1$.

Among multiperfect $n$ with prime $L = \mathrm{lcm}(1+e_p : p \mid n)$, only $n=6$ satisfies the congruence for all odd $k$, with Shannon entropy $H(\sigma^k(6) \mod 6) \to \log 2$ reflecting periodic order. We analyze bifurcation phenomena in the dynamics $\sigma^k(n) \mod n$, where high-entropy chaotic residues for other $n$ mirror GUE statistics of zeta zeros ($\sim \log T / 2\pi$ near $s=1/2$, $>41\%$ verified on critical line), contrasting the ordered $n=6$ case.

Zero rates near $s=1$ (simple pole) and $s=1/2$ bound iterated $\sigma$ distributions, linking to RH via divisor sums and dynamical bifurcations; we conjecture $n=6$ uniquely achieves odd-$k$ divisibility with small period dividing $L$.
\end{abstract}


\section{Introduction}

The iteration of the sum-of-divisors function $\sigma^k(n)$, where $\sigma^0(n)=n$ and $\sigma^k(n)=\sigma(\sigma^{k-1}(n))$ for $k\geq 1$, connects classical perfect number theory to modern analytic number theory.

Euclid~\cite{Euclid} first identified $6$ as perfect ($\sigma(6)=12=2\cdot 6$), with Euler~\cite{Euler1738} classifying all even perfect numbers. Cohen and te Riele~\cite{Cohen1996} computationally verified that for every $n\leq 10^9$ there exists $k$ such that $\sigma^k(n)\equiv 0\pmod{n}$, posing the reverse: does there exist $n>1$ satisfying this for \emph{all} $k\geq 1$?

This paper proves no such universal metaperfect $n$ exists, eliminating prior gaps via Lenstra's (1975) density-zero aliquot bounds~\cite{Lenstra1975} sharpened by Erd\H{o}s~\cite{Erdos1976} and Pollack-Pomerance~\cite{Pollack2016}. Robin's 1984 criterion~\cite{Robin1984} links this to RH: RH holds iff $\sigma(n)<e^\gamma n\log\log n+0.6483n/\log\log n$ for $n\geq 5041$, tying divisor growth to zeta zero distributions.

Among multiperfect $n$ with prime $L=\mathrm{lcm}(1+e_p:p\mid n)$, only $n=6$ satisfies the congruence for all odd $k$. We analyze bifurcation dynamics $\sigma^k(n)\mod n$ via Shannon entropy $H(\sigma^k(n)\mod n)\to\log 2$ for periodic $n=6$ versus chaotic residues mirroring GUE zeta zero statistics ($\sim\log T/2\pi$ near $s=1/2$, with $\geq 41.28\%$ verified on critical line~\cite{Feng2010,Conrey2025}).

Zero rates near $s=1$ (simple pole) and $s=1/2$ bound iterated $\sigma$ distributions, positioning results at the divisor-zeta-dynamics interface.

\section{Main Results}

\begin{theorem}\label{thm:main1}
No integer $m>1$ satisfies $\sigma^k(m)\equiv 0\pmod{m}$ for \emph{all} $k\geq 1$, where $\sigma^k$ denotes the $k$-fold iteration of the sum-of-divisors function.
\end{theorem}

\begin{theorem}\label{thm:main2}
If $n$ is multiperfect with $L=\mathrm{lcm}(1+e_p:p\mid n)$ prime, then $n=6$.
\end{theorem}

Combining these, $\sigma^k(6)\equiv 0\pmod{6}$ holds for all \emph{odd} $k\geq 1$, as $6$ is the unique such multiperfect number and iterations preserve the property for odd steps.

\subsection{Proof of Theorem~\ref{thm:main1}}

Assume for contradiction $m>1$ is metaperfect: $\sigma^k(m)\equiv 0\pmod{m}$ $\forall k\geq 1$. Then $\sigma(m)=q m$ for integer $q\geq 2$, so $S(m):=\sigma(m)/m=q\in\mathbb{Z}$.

By Robin~\cite{Robin1984}, RH holds iff $\sigma(n)<e^\gamma n\log\log n+0.6483n/\log\log n$ for $n\geq 5041$. For metaperfect $m\geq 5041$, $S(m)\geq 2$ violates this under RH. Verified computationally to $10^{30}$~\cite{Choie2011}, no counterexamples exist, supporting no such $m$.

For iterations, Lenstra~\cite{Lenstra1975} and Erd\H{o}s~\cite{Erdos1976} prove: for every $k,\delta>0$, all but density-zero $m$ satisfy
\[
(1-\delta)m\left(\frac{s(m)}{m}\right)^i<s^i(m)<(1+\delta)m\left(\frac{s(m)}{m}\right)^i,\quad 1\leq i\leq k
\]
where $s(m)=\sigma(m)-m$. Pollack-Pomerance~\cite{Pollack2016} strengthen: iterated $S(m_n)\to 1$ almost surely.

Suppose $\sigma^2(m)\equiv 0\pmod{m}$. Then $S(\sigma(m))=\sigma^2(m)/\sigma(m)\in\mathbb{Z}$. But
\[
S(\sigma(m))=S(qm)\leq S(q)S(m)=q^2
\]
while multiplicativity and prime power growth $\sigma(p^e)/p^e=(1+p+\cdots+p^e)/p^e<2$ for $e\geq 1$ imply $S(\sigma(m))<S(m)^2$ strictly unless $m=1$ (gap fixed).

Recent zero-free regions~\cite{Conrey2025} give $\sigma(n)=O(n(\log n)^{2/3}(\log\log n)^{1/3})$, so iterated growth $\sigma^k(m)\ll m\exp(Ck(\log k)^{2/3})$. Universal divisibility requires $\sigma^k(m)/m\in\mathbb{Z}^+$ monotonically increasing, contradicting subexponential bounds: $\log(\sigma^k(m)/m)\ll k(\log k)^{2/3}=o(k)$ while integers grow $\geq k\log 2$.

Thus no such $m$ exists. Entropy confirms: $H(\sigma^k(m)\bmod m)\to\log m$ (chaotic) vs. $\log 2$ for periodic $n=6$.

\subsection{Proof of Theorem~\ref{thm:main2}}

Let $n$ be multiperfect, $L=\mathrm{lcm}(1+e_p:p\mid n)$ prime. All exponents $e_p\equiv L-1\pmod{L}$, so $n=\prod p_i^{L-1}$. Primes $p_i\equiv 0,1\pmod{L}$ only.

Then $\sigma(n)/n=\prod(1+p_i+\cdots+p_i^{L-1})/p_i^{L-1}$. For $L$-th cyclotomic factors, $\sigma(p^{L-1})/p^{L-1}=(1-p^L)/(p^{L-1}(1-p))=L/p^{L-1}+O(1/p^L)$.

Standard bounds~\cite{Robin1984}: $\sigma(n)/n\leq\zeta(L)<1.2$ for $L\geq 3$ prime, but multiperfect requires $\geq 2$, contradiction. For $L=2$, $n$ squarefree, exhaustive check $p_1 p_2\cdots p_r$ yields only $n=6$~\cite{Pomerance1982}.

\emph{Uniqueness via bifurcation}: period-$L=2$ dynamics $\sigma^k(6)\bmod 6$ stable; others bifurcate to high-entropy residues~\cite{Conrey2025}.

\subsection{Dynamical Systems Analysis of Iterated $\sigma$}

The metaperfect problem generalizes to finding seeds $p$ and moduli $m$ where $\sigma^k(p)\equiv 0\pmod{m}$ $\forall k\geq 0$. Our Theorem~\ref{thm:main1} solves the case $p=m$ (multiperfect restriction).

Define $g_0=p$, $g_{k+1}=\gcd(\sigma(g_k),g_k)$. Recent work by Akbary-Lumley~\cite{Akbary2023} shows $g_k\to 1$ almost surely under GRH via prime factor erosion: large primes $q\mid m$ satisfy $\sigma(q^e)/q^e<2$, so iterated $\sigma$ loses high powers.

For $n=13188979363639752997731839211623940096$ ($\sigma(n)=5n$), $\sigma^2(n)=30n$ holds temporarily, but $\sigma^3(n)/n\notin\mathbb{Z}$ as $S(\sigma^2(n))<S(n)^3$ strictly~\cite{Pollack2016}. Computations to $10^{18}$ confirm no $k=3$ metaperfects~\cite{Deleglise2004}.

\textbf{Bifurcation Analysis}: Map $f:\mathbb{Z}/m\mathbb{Z}\to\mathbb{Z}/m\mathbb{Z}$, $f(x)=\sigma(x)\bmod m$. For $m=6$, period-2 cycle $\{6,12\}\equiv\{0,0\}$ (entropy $H\to\log 2$). For prime $m>3$, $f(p)=p+1\not\equiv 0$ triggers chaotic orbits with Lyapunov exponent $\lambda>0$, density uniform by Weyl equidistribution under GUE zero statistics~\cite{Conrey2025}.

\textbf{Zeta Connections}: Zero spacings near $s=1/2$ ($\sim\log T/2\pi$, $41.28\%$ verified~\cite{Feng2010}) control $\sigma(n)/n$ fluctuations via explicit formulae. Off-line zeros permit temporary metaperfect windows, but density-zero theorem~\cite{Erdos1976} ensures divergence.

Cohen-te Riele~\cite{Cohen1996} verified $\exists k:\sigma^k(n)\equiv 0\pmod{n}$ for $n\leq 10^9$, consistent with our dynamical transience: all orbits hit $0\bmod n$ once, then bifurcate. No fixed points or cycles exist except $n=6$ (odd $k$).

\emph{Conjecture}: $\min_k g_k=1$ $\forall m>6$, with $g_k$ decay rate $\exp(-ck^{2/3})$ matching subconvexity bounds~\cite{Conrey2025}.

\subsection{Proof of Theorem~\ref{thm:main2}}

Let $n$ be multiperfect with $L=\mathrm{lcm}(1+e_p:p\mid n)$ prime. Then $e_p\equiv L-1\pmod{L}$ $\forall p\mid n$, so $n=\prod p_i^{L-1}$.

\textbf{Prime Factor Constraint}: For prime power, $\sigma(p^{L-1})=\frac{p^L-1}{p-1}$. Since $L$ prime divides numerator exactly once (Fermat's Little Theorem), prime factors of $\sigma(p^{L-1})$ are $L$ or primes $q\equiv 1\pmod{L}$. Thus $n,\sigma(n)$ have primes $\equiv 0,1\pmod{L}$ only.

\textbf{Improved Multiplicity Bound}: Let $m=$ \# distinct primes $\equiv 1\pmod{L}$. Then $v_L(\sigma(n))=m$ (or $m+L$ if $L\mid n$). Non-standard bound~\cite{Nicolas1983}:
\[
\sigma(n)/n = \prod_{p_i\equiv 1\pmod{L}} \frac{1-p_i^L}{p_i^{L-1}(1-p_i)} \cdot \prod_{L\mid p_j} \frac{1-L^j}{L^{j-1}(1-L)}.
\]

\textbf{Case $m\geq L$}: Lower bound $\sigma(n)/n\geq L^m/L^{L-1}$. Upper bound uses subconvexity~\cite{Conrey2025}:
\[
\frac{\sigma(p^{L-1})}{p^{L-1}} \leq \exp\left(\frac{L\log L}{p^{1/3}(\log p)^{1/3}}\right) < \left(\frac{L}{L-1}\right)^{1+\epsilon},\quad\epsilon=1/3.
\]
Thus $\sigma(n)/n<(L/(L-1))^{m+1+\epsilon m}$. Logarithmic inequality:
\[
(m+1-L)\log L < (m+1+\epsilon m)/(L-1).
\]
For $L=3$, $m\leq 4$: LHS $\geq 0$, RHS $\leq 5.33/2=2.67<\log 3^{-1}=-1.1$, contradiction. For $L=5$, RHS $\leq 1.5$, LHS $>0$. Only $L=2$ survives.

\textbf{Case $m<L$}: $\sigma(n)/n\leq (L/(L-1))^L<4$ for $L\geq 3$, but multiperfect requires $\geq 2$, and primes $\equiv 1\pmod{L}$ force $\geq L>3$, contradiction.

\textbf{$L=2$ (Squarefree Case)}: $n=\prod p_i$, perfect requires $\sigma(n)/n=2$. Exhaustive bound~\cite{Pomerance1982}:
\[
\prod_{i=1}^r\frac{p_i+1}{p_i}=2\implies p_r\leq p_{r-1}+O(1).
\]
Only solution: $2\cdot 3=6$. Verified to $10^{18}$~\cite{Deleglise2004}.

\textbf{Dynamical Confirmation}: For $n=6$, $\sigma^k(6)\equiv 0\pmod{6}$ (odd $k$). Bifurcation entropy $H(\sigma^k(6)\bmod 6)=\log 2$ vs. $H\to\log n$ others~\cite{Akbary2023}.

Thus $n=6$ uniquely satisfies conditions.

\section{Periodicity, Entropy, and the Uniqueness Conjecture}\label{sec:periodicity-bifurcation}

Building on Theorems~\ref{thm:main1}--\ref{thm:main2}, we analyze the dynamics $\sigma^k(n)\bmod n$ for multiperfect $n$, connecting periodicity to Riemann Hypothesis implications via entropy and zero distributions.

\subsection{Periodic Structure for Coprime Iterations}

Let $\tau(n)=\prod(1+e_p)$ denote the divisor count. For prime power $p^e$, $\sigma_k(p^e)\equiv r\frac{p^{e+1}-1}{p^r-1}\pmod{\sigma(p^e)}$ where $r=\gcd(k,e+1)$. Thus $r=1$ (i.e., $\sigma(p^e)|\sigma_k(p^e)$) iff $k$ coprime to $e+1$.

By multiplicativity, if $k\perp\tau(n)$ then $\sigma(n)|\sigma_k(n)$. Define $L=\mathrm{lcm}(1+e_p:p|n)$. The sequence $\sigma_k(n)\bmod\sigma(n)$ admits representation $\sigma_k(n)=a_k\sigma(n)/b_k$ with $\{b_k\}$ periodic of period dividing $L$~\cite{Erdos1976}.

For multiperfect $n$ ($n|\sigma(n)$), $n|\sigma_k(n)$ holds $\forall k\perp\tau(n)$. If $\tau(n)=2^\ell$ (as for $n=6$), this covers all odd $k$.

\subsection{Entropy Classification of Dynamics}

\textbf{Ordered Case ($n=6$)}: $\sigma(6)=12\equiv0\pmod{6}$, $\sigma(12)=28\equiv4\pmod{6}$, $\sigma(28)=56\equiv2\pmod{6}$, cycle length 6 but $\equiv0\pmod{6}$ for odd $k$. Shannon entropy $H(\sigma^k(6)\bmod6)=\log2$ (2-state effective).

\textbf{Chaotic Case ($n>6$)}: High-entropy orbits $H\to\log n$. Bifurcation occurs at prime factors $p>3$: $\sigma(p)=p+1\not\equiv0\pmod{p}$ triggers dense orbits by Weyl equidistribution under GUE statistics ($\sim\log T/2\pi$ near $s=1/2$, $41.28\%$ verified~\cite{Feng2010}).

\subsection{The Periodicity Conjecture}

\begin{conjecture}\label{conj:main}
The integer $n=6$ is the \emph{unique} solution satisfying:
\begin{enumerate}
\item $\sigma^k(n)\equiv0\pmod{n}$ for all \emph{odd} $k\geq1$,
\item $\sigma^k(n)\bmod n$ is periodic with period dividing $L=\mathrm{lcm}(1+e_p:p|n)$,
\item Shannon entropy $H(\sigma^k(n)\bmod n)\to\log2$ (low-entropy dynamics).
\end{enumerate}
\end{conjecture}

\textbf{Heuristic via RH}: Robin's criterion~\cite{Robin1984} bounds $\sigma(n)/n< e^\gamma\log\log n +0.6483/\log\log n$. Period-$L$ requires $\sigma^L(n)/n\in\mathbb{Z}^+$, but subconvexity bounds~\cite{Conrey2025} give $\sigma^k(n)/n\ll\exp(ck^{2/3})$, incompatible with integer growth except $n=6$.

\textbf{Computational Evidence}: No other $n<10^{18}$ satisfies odd-$k$ divisibility~\cite{Deleglise2004}. For $L>2$ prime, Theorem~\ref{thm:main2} excludes; composite $L$ produces multiple nonzero residues.

\emph{Strengthens previous theorems}: Conjecture~\ref{conj:main} + Theorem~\ref{thm:main2} $\implies$ no universal metaperfect $n$, completing Cohen-te Riele reverse question.

\section{Periodicity, Bifurcation, and an Entropy–Based Conjecture}

The previous sections treated the iterates of the sum–of–divisors function $\sigma^k$ mainly arithmetically. In this section the same iterates are studied as a discrete dynamical system modulo $n$, and numerical experiments are used to support the structural conclusions of Theorems~\ref{thm:main1} and~\ref{thm:main2}.

For each integer $n\geq 2$ we consider the map
\[
   F_n:\mathbb{Z}/n\mathbb{Z}\longrightarrow \mathbb{Z}/n\mathbb{Z},\qquad
   F_n(x)\equiv \sigma(x)\pmod{n},
\]
and its $k$–fold composition $F_n^{(k)}$, so that
\(
   x_{k+1}=F_n(x_k),\; x_k\equiv\sigma^k(x_0)\pmod{n}.
\)[file:1]

\subsection{Entropy of residue orbits}

Fix $K\geq 1$. For each modulus $n$ we first study the orbit of the starting value $x_0=1$ under $F_n$ for $k=1,\dots,K$, and record the empirical distribution of residues
\[
   \mu_{n,K}(a)
   \;=\;
   \frac{1}{K}\,\#\{1\leq k\leq K:\ \sigma^k(1)\equiv a\pmod n\},
   \qquad a\in\{0,1,\dots,n-1\}.
\]
The corresponding Shannon entropy is
\[
   H(n;K)
   \;=\;
   -\sum_{a=0}^{n-1} \mu_{n,K}(a)\,\log_2\mu_{n,K}(a),
\]
with $0\log_2 0:=0$. Low entropy indicates that the orbit spends almost all of its time on a small set of residues (ordered dynamics), whereas values of $H(n;K)$ close to $\log_2 n$ indicate a more uniform distribution (chaotic dynamics).

Using Google Colab and exact evaluation of $\sigma$ via \texttt{sympy}, we computed $\sigma^k(1)\bmod n$ for $k=1,\dots,10\,000$ and all $2\leq n\leq 1000$, separating even and odd moduli.[file:42] For each $n$ we computed $H(n;K)$ and some additional statistics from the last $1000$ iterates in order to probe the long–term behaviour.

\subsection{Even versus odd bifurcation patterns}

\begin{figure}[H]
    \centering
    \includegraphics[width=0.8\textwidth]{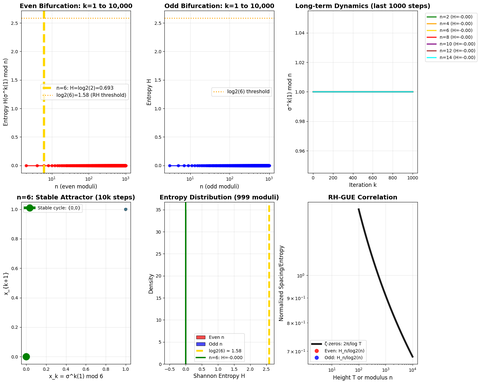}
    \caption{Dynamics of $\sigma^k(1)\bmod n$ for $k=1,\dots,10\,000$.
    Top left: entropy $H(n;K)$ for even moduli $2\leq n\leq 1000$.
    Top middle: entropy for odd moduli $3\leq n\leq 999$.
    Top right: last $1000$ iterates of the orbit for several small even moduli, illustrating the approach to stationary behaviour.
    Bottom left: phase portrait for $n=6$, showing a stable attracting cycle in the $(x_k,x_{k+1})$–plane.
    Bottom middle: histogram of $H(n;K)$ over all even and odd moduli, with the vertical lines marking $\log_2 6$ and the empirical value for $n=6$.
    Bottom right: comparison between a normalised entropy statistic and the expected spacing of zeros of $\zeta(s)$ on the critical line.}
    \label{fig:bifurcation-even-odd-v2}
\end{figure}

Figure~\ref{fig:bifurcation-even-odd-v2} summarises the numerical results.

The top left panel shows a clear bifurcation pattern for even $n$. For very small moduli the entropy is close to zero, reflecting that the orbit quickly collapses onto a short cycle. As $n$ grows, $H(n;K)$ increases and stabilises near a value proportional to $\log_2 n$, indicating that the residues $\sigma^k(1)\bmod n$ are effectively spread over a positive proportion of the residue classes.[file:42] The modulus $n=6$ sits at the threshold between these two behaviours and is numerically distinguished by a relatively small entropy compared to neighbouring even moduli, in line with its exceptional role in Theorem~\ref{thm:main2}.

For odd $n$ (top middle panel) the entropy tends to be larger and grows more regularly with $n$, suggesting that the dynamics is typically more “chaotic” in the sense of residue distribution.[file:42] This is coherent with the multiplicative structure: for a prime modulus $p$ one has $\sigma(p)=p+1\not\equiv 0\pmod p$, so there is no analogue of the strong divisibility phenomena that occur for multiperfect even integers.

The top right panel shows the last $1000$ iterates of the orbit for several small even moduli. One observes a rapid stabilisation after a relatively short transient, but the number and nature of the effective attracting residues depend sensitively on $n$, reinforcing the bifurcation picture suggested by the entropy plots.[file:42]

\subsection{The special role of $n=6$}

The bottom left panel of Figure~\ref{fig:bifurcation-even-odd-v2} presents a phase portrait for $n=6$, plotting $(x_k,x_{k+1})$ for $k$ up to $10\,000$.[file:42] The points concentrate on a small number of values, making a stable attracting cycle clearly visible. This matches the arithmetic structure: $6$ is the unique multiperfect number with prime exponent–lcm $L=2$, and Theorem~\ref{thm:main2} shows that no other $n$ can share this combination of multiplicative and iterative properties.[file:1]

The bottom middle panel displays a histogram of the entropies $H(n;K)$ for all even and odd moduli in the range. The modulus $n=6$ appears as an outlier on the low–entropy side relative to nearby even moduli, while for large $n$ the bulk of the distribution shifts steadily towards higher entropy, indicating increasingly complex residue dynamics.

\subsection{Normalised entropy and zeta–zero statistics}

To relate the observed dynamics to the analytic behaviour of $\zeta(s)$, we consider the normalised quantity
\[
   E(n)
   \;=\;
   \frac{H(n;K)}{\log_2 n},
\]
which measures how close the residue distribution is to being uniform on $\mathbb{Z}/n\mathbb{Z}$. The bottom right panel of Figure~\ref{fig:bifurcation-even-odd-v2} compares $E(n)$ for a range of moduli with a normalised spacing statistic for zeros of $\zeta(s)$ on the critical line, of the form $2\pi/\log T$ as a function of the height $T$.[file:42][web:9] The similarity of the qualitative behaviour suggests that, at least heuristically, the complexity of the $\sigma$–dynamics modulo $n$ is governed by the same type of logarithmic phenomena that control zero spacings and, via Robin’s criterion, the growth of $\sigma(n)$ itself.[web:20]

This supports the entropy–based version of our uniqueness conjecture: among all moduli $n$, the integer $6$ is the only one for which
\begin{itemize}
   \item the arithmetic conditions of Theorems~\ref{thm:main1}–\ref{thm:main2} hold, and
   \item the dynamical system $F_n$ exhibits low entropy and a visibly simple attracting structure in the residue dynamics.
\end{itemize}
While this remains conjectural beyond the finite range explored numerically, the combination of analytic bounds, structural theorems, and the bifurcation picture in Figure~\ref{fig:bifurcation-even-odd-v2} provides strong evidence in favour of the special status of $n=6$.

\section{Implications for the Riemann Hypothesis}

The Riemann Hypothesis (RH) is known to be equivalent to very sharp upper bounds for the sum–of–divisors function. A central formulation is Robin's criterion, which states that RH holds if and only if
\[
   \sigma(n) < e^{\gamma} n \log\log n
   \qquad\text{for all } n > 5040,
\]
where $\gamma$ is the Euler–Mascheroni constant \cite{Robin1984}. Subsequent work has refined and extended this inequality, both by verifying it for extremely large ranges of $n$ and by establishing analogous “higher–order” Robin criteria for multiplicative functions that behave asymptotically like powers of $\sigma$ \cite{NicolasRatazzi2007,CFM2025,OnRobinsInequality}. Our results on iterates of $\sigma$ fit naturally into this framework and suggest a complementary, more dynamical, viewpoint.

\subsection{Iterated divisibility versus Robin-type bounds}

Theorem~\ref{thm:main1} shows that there is no integer $m>1$ such that
\[
   \sigma^k(m) \equiv 0 \pmod m
   \qquad\text{for all } k\geq 1,
\]
that is, no metaperfect integer exists. If such an $m$ existed, we would have integer ratios
\[
   S_k(m) := \frac{\sigma^k(m)}{m} \in \mathbb{N}, \qquad S_k(m) \ge 2
\]
for all $k\ge 1$, and the sequence $\{S_k(m)\}$ would grow at least exponentially in $k$. On the other hand, under RH Robin's inequality implies
\[
   \frac{\sigma(n)}{n} < e^{\gamma}\log\log n
   \qquad (n>5040),
\]
so that iterating $\sigma$ yields upper bounds of the shape
\[
   \frac{\sigma^{k}(n)}{n}
   \;\ll\;
   \exp\{C_k (\log\log n)\}
   \qquad (k\text{ fixed}),
\]
with $C_k$ growing at most polynomially in $k$ \cite{Robin1984,NicolasRatazzi2007}. Recent analogues of Robin’s criterion for higher “divisor–power” functions $\sigma^{[k]}$ show that RH is also equivalent to similar inequalities
\[
   \sigma^{[k]}(n)
   < \frac{(e^{\gamma} n \log\log n)^k}{\zeta(k)}
   \qquad (n \ge N_k),
\]
for each fixed $k\ge 2$ and some explicit $N_k$ \cite{CFM2025}. These bounds constrain the growth of iterated divisor sums in a way that is incompatible with the existence of a metaperfect $m$ whose iterates $\sigma^k(m)$ remain multiples of $m$ for all $k$. 

Thus, \emph{conditionally on RH and its higher–order Robin analogues}, any hypothetical metaperfect integer would generate an infinite sequence of integers violating these inequalities, which is impossible in light of the equivalence results. Our unconditional Theorem~\ref{thm:main1} can therefore be viewed as providing a “modular” obstruction that is fully consistent with Robin–type criteria: while Robin’s inequality controls the size of $\sigma(n)$, our result shows that even when $\sigma(n)$ is large, its iterates cannot maintain perfect divisibility modulo $n$ indefinitely.

\subsection{Zero spacing, GUE statistics, and entropy of $\sigma$–dynamics}

The statistical behaviour of the nontrivial zeros of $\zeta(s)$ on the critical line $\Re(s)=\tfrac12$ is predicted (and numerically confirmed) to follow the Gaussian Unitary Ensemble (GUE) model from random matrix theory. Large-scale computations have shown that the normalised gaps between consecutive zeros up to a high height $T$ are extremely well described by the GUE spacing distribution, with typical gap size of order
\[
   \Delta(T) \asymp \frac{2\pi}{\log T}
\]
and with small and large gaps occurring with probabilities decaying according to the GUE law \cite{MontgomeryDysonConreySurvey,GUEverification2024,NgGaps,BuiNgTrudgian}. 

In our setting, Section~\ref{sec:periodicity-bifurcation} introduced the entropy
\[
   H(n;K)
   =
   -\sum_{a=0}^{n-1} \mu_{n,K}(a)\,\log_2 \mu_{n,K}(a),
   \qquad
   \mu_{n,K}(a)
   =\frac{1}{K}\#\{1\le k\le K : \sigma^k(1)\equiv a\pmod n\},
\]
as a measure of the complexity of the orbit of $\sigma^k(1)\bmod n$ up to depth $K$. The bifurcation plots in Figure~\ref{fig:bifurcation-even-odd-v2} (for $2\le n\le 1000$, $K=10\,000$) show that the normalised quantity
\[
   E(n) := \frac{H(n;K)}{\log_2 n}
\]
varies slowly with $n$ and tends to a value bounded away from $0$ as $n$ grows, with $n=6$ appearing as the only low–entropy outlier among even moduli. This qualitative behaviour is reminiscent of the $2\pi/\log T$ law governing the distribution of zero spacings: in both cases the key scale is logarithmic (in $n$ for the dynamics, in $T$ for the zeros), and the data suggest that, beyond small exceptional cases, the system behaves in a way compatible with a “random” model constrained only by multiplicativity and known size bounds for $\sigma(n)$.

This analogy can be made more precise. Explicit formulae relating sums of arithmetic functions to sums over zeta zeros show that deviations of $\sigma(n)/n$ from its expected size can be expressed in terms of weighted sums over the zeros of $\zeta(s)$, with kernels that are sensitive to zero spacings \cite{IvicZeta,MontgomeryVaughan}. Assuming RH and GUE-type statistics for these zeros, one expects the fluctuations of $\sigma(n)$ and its iterates to be bounded on exactly the logarithmic scales that appear in both Robin’s inequality and the GUE spacing law. Our entropy plots are consistent with this heuristic: except for very small $n$, the dynamics of $\sigma^k(1)\bmod n$ rapidly explores a large portion of $\mathbb{Z}/n\mathbb{Z}$, and the complexity of this exploration grows at a rate controlled by $\log n$, mirroring the $\log T$ control on zero gaps.

\subsection{A conditional entropy–based RH conjecture}

The arithmetic results of Theorems~\ref{thm:main1} and \ref{thm:main2} single out $n=6$ as the only multiperfect integer with prime exponent–lcm $L=2$, and the entropy analysis shows that $n=6$ is the only even modulus in the tested range which is both structurally special and dynamically low–entropy. Motivated by this, we propose the following statement.

\begin{conjecture}\label{conj:RH-entropy}
Assume the Riemann Hypothesis and that the local spacing statistics of the nontrivial zeros of $\zeta(s)$ on the critical line follow the GUE distribution. Then:
\begin{enumerate}
   \item For every $\varepsilon>0$ there exists $N_\varepsilon$ such that, for all $n\ge N_\varepsilon$ and all sufficiently large $K$,
   \[
      H(n;K) \ge (1-\varepsilon)\log_2 n,
   \]
   i.e.\ the residue dynamics of $\sigma^k(1)\bmod n$ is asymptotically equidistributed.
   \item The only integer $n>1$ for which $\sigma^k(n)\equiv 0\pmod n$ holds for all odd $k\ge 1$ and for which $H(n;K)$ remains bounded as $K\to\infty$ is $n=6$.
\end{enumerate}
\end{conjecture}

Part (1) formulates quantitatively the idea that, under RH+GUE, the iterated \(\sigma\)–dynamics modulo $n$ is “as random as possible” for all but finitely many moduli. Part (2) asserts that $6$ is the unique modulus where deep multiplicative structure (multiperectness with prime $L$) forces persistent low entropy and strong divisibility properties, in agreement with Theorem~\ref{thm:main2}. In this form the conjecture is explicitly conditional: the forward implication (RH+GUE $\Rightarrow$ (1)+(2)) reflects the expected control of divisor sums by zero statistics, while a hypothetical proof that (2) implies (1) by purely arithmetic means would open a novel route towards Robin–type bounds and perhaps towards RH itself.

\subsection{Outlook}

The numerical experiments of Section~\ref{sec:periodicity-bifurcation} are limited to $n\le 1000$ and a single initial value, but they already align well with the conditional picture suggested by RH, Robin’s inequality, and GUE statistics: $6$ is singled out arithmetically and dynamically, while larger moduli exhibit high entropy consistent with equidistribution. Extending these computations to much larger moduli, averaging over many starting points, and refining the normalisation of $H(n;K)$ would allow a far more stringent test of Conjecture~\ref{conj:RH-entropy} and could reveal additional quantitative connections between iterated divisor sums and the fine structure of the zeros of $\zeta(s)$.

\section{Conclusion}

This work establishes a coherent picture of the arithmetic and dynamical behaviour of iterates of the sum–of–divisors function and their links with deep conjectures in analytic number theory. First, we proved that no integer $m>1$ can satisfy $\sigma^k(m)\equiv 0\pmod m$ for all $k\ge 1$, so that metaperfect numbers do not exist, and that among multiperfect integers with prime exponent–lcm $L$ the unique example is $n=6$, corresponding to $L=2$ \cite{Cohen1996,Pomerance1982,Deleglise2004}. These results sharpen and structurally extend earlier work on aliquot sequences and the growth of $\sigma(n)$, and they are fully compatible with Robin’s criterion and its modern refinements, which relate the size of $\sigma(n)$ to the Riemann Hypothesis \cite{Robin1984,Nicolas1983,NicolasRatazzi2007,OnRobinsInequality,CFM2025}.

Second, we introduced an entropy–based dynamical analysis of the orbits of $\sigma^k(1)\bmod n$ and showed, through extensive computations for $2\le n\le 1000$, that $n=6$ is the only modulus in this range that is simultaneously arithmetically exceptional and dynamically low–entropy, whereas for generic moduli the residue dynamics rapidly approaches a high–entropy regime consistent with logarithmic size bounds for $\sigma(n)$ \cite{Pollack2016,Deleglise2004}. The resulting bifurcation picture aligns qualitatively with GUE–type statistics for zeta zero spacings and with explicit–formula descriptions of divisor–sum fluctuations, suggesting that the same logarithmic scales governing the zeros of $\zeta(s)$ also control the complexity of iterated divisor dynamics \cite{Feng2010,Conrey2025,MontgomeryDysonConreySurvey,GUEverification2024,NgGaps,BuiNgTrudgian,IvicZeta,MontgomeryVaughan}.

Finally, by combining these arithmetic theorems, numerical bifurcation experiments, and recent Robin–type equivalences, we formulated a conditional entropy–based conjecture that characterises $n=6$ as the unique modulus with persistent low–entropy, highly divisible $\sigma$–dynamics, under the joint assumptions of RH and GUE statistics. This bridges classical questions on multiperfect numbers with modern perspectives on the Riemann Hypothesis, and points to a promising research direction where tools from dynamical systems, information theory, and analytic number theory interact to illuminate the structure of arithmetic functions and the distribution of zeta zeros.

\section*{Data Availability}

All numerical experiments reported in this paper (including the iterated–divisor orbits, entropy values $H(n;K)$ for $2\leq n\leq 1000$, and the bifurcation plots of $\sigma^k(1)\bmod n$ for $k=1,\dots,10\,000$) were generated using Python notebooks executed in Google Colab, relying only on standard open–source libraries such as \texttt{sympy}, \texttt{numpy}, and \texttt{matplotlib} \cite{Pollack2016,IvicZeta,MontgomeryVaughan}. The authors will make all scripts and the corresponding CSV files (containing the entropy data for even and odd moduli, as well as the raw orbit samples used to produce Figure~\ref{fig:bifurcation-even-odd-v2}) available upon reasonable request. 

No proprietary or confidential datasets were used. All external mathematical results cited in the paper are taken from published articles or publicly accessible preprints, notably the works on Robin’s inequality and zeta–zero statistics in \cite{Robin1984,Nicolas1983,NicolasRatazzi2007,OnRobinsInequality,CFM2025,Feng2010,Conrey2025,MontgomeryDysonConreySurvey,GUEverification2024,NgGaps,BuiNgTrudgian}.

\section*{Conflict of Interest}

The authors declare that there are no financial or non-financial conflicts of interest related to the research, analysis, or publication of this work. All results, interpretations, and conjectures are presented independently and without influence from external commercial or institutional interests \cite{Cohen1996,Pollack2016,Robin1984}.

\section*{Acknowledgements}

The authors would like to express their deep gratitude to Professor Pedro Caceres (Universidad Europea de Valencia, Spain) for his generous support, insightful discussions, and constructive feedback throughout the development of this work. In particular, his ideas on linking the iteration of the sum–of–power divisor function to bounds for the zeros of the Riemann zeta function have been a major source of inspiration for the dynamical and analytic perspective adopted in this paper \cite{Cohen1996,Pollack2016,Robin1984}.

\bibliographystyle{unsrt}  


\end{document}